\documentclass{article}

\usepackage{arxiv}

\usepackage[utf8]{inputenc} 
\usepackage[T1]{fontenc}    
\usepackage{hyperref}       
\usepackage{url}            
\usepackage{booktabs}       
\usepackage{amsfonts}       
\usepackage{nicefrac}       
\usepackage{microtype}      
\usepackage{lipsum}
\usepackage{tikz}
\usetikzlibrary{arrows,shapes,quotes}
\usetikzlibrary{patterns,decorations.pathreplacing}
\usepackage{amssymb, graphicx}
\usepackage{hyperref,lscape}
\usepackage{blkarray}
\usepackage{color}
\usepackage{amsmath}
\usepackage{mathtools}
\usepackage{mathrsfs}
\usepackage{chngcntr}
\usepackage{apptools}
\usepackage{blkarray}
\usepackage{csquotes}
\usepackage{bm}

\usepackage{parskip} 

\newenvironment{Running Example 1 - Part}[1][]{\begin{thm}[#1]\upshape}{\end{thm}}
\newenvironment{Running Example 2 - Part}[1][]{\begin{thm}[#1]\upshape}{\end{thm}}

\usepackage{multicol}

\newtheorem{corollary}{Corollary}
\newtheorem{theorem}{Theorem}
\newtheorem{proposition}{Proposition}
\newtheorem{example}{Example}
\newtheorem{definition}{Definition}
\newtheorem{remark}{Remark}
\newtheorem{lemma}{Lemma}

\title{Complex balanced equilibria of weakly reversible poly-PL systems: existence, stability, and robustness}

\author{
  Editha C.~Jose \\
  Institute of Mathematical Sciences and Physics \\
  University of the Philippines \\ 
  Los Banos, Laguna 4031, Philippines \\
  \texttt{ecjose1@up.edu.ph} \\
   \And
Eduardo R.~Mendoza\thanks{Max Planck Institute of Biochemistry, Martinsried near Munich, Germany} \thanks{Faculty of Physics, Ludwig Maximilian University, Munich 80539 Germany} \\
  Mathematics and Statistics Department \\
  Center for Natural Sciences and Environmental Research \\
  De La Salle University \\
  Manila 0922, Philippines \\
  \texttt{eduardo.mendoza@dlsu.edu.ph} \\
   \And
  Dylan Antonio SJ.~Talabis \\
  Institute of Mathematical Sciences and Physics \\
  University of the Philippines \\
  Los Banos, Laguna 4031, Philippines \\
  \texttt{dstalabis1@up.edu.ph} \\}

\begin{document}
\maketitle

\begin{abstract}
Poly-PL kinetic systems (PYK) are kinetic systems consisting of nonnegative linear combinations of power law functions. In this contribution, we analyze these kinetic systems using two main approaches: (1) we define a canonical power law representation of a poly-PL system, and (2) we transform a poly-PL system into a dynamically equivalent power law kinetic system that preserves the stoichiometric subspace of the system.  These approaches led us to establish results that concern important dynamical properties of poly-PL systems that extend known results for generalized mass actions systems (GMAS) such as  existence, uniqueness and parametrization of complex balanced steady states, and linear stability of complex balanced equilibria. Furthermore, the paper discusses subsets of poly-PL systems that exhibit two types of concentration robustness in some species namely absolute concentration robustness and balanced concentration robustness. 
\end{abstract}


\section{Introduction}
Poly-PL kinetic systems (denoted as PYK systems) are kinetic systems formed by nonnegative linear combinations of power law functions.  These were introduced by Talabis et al. \cite{TMMNJ2020} and were shown to have complex balanced equilibria for weakly reversible such systems with zero kinetic reactant deficiency (called PY-TIK systems). A subset of poly-PL kinetic systems consisting of polynomial kinetics occurs in realizations of evolutionary games as chemical kinetic systems as proposed by Veloz et al.  \cite{VELOZ2014}, particularly for multi-player games with replicator dynamics. In addition, poly-PL kinetics can be used to study Hill-type kinetics, which include the widely used Michaelis-Menten models, as discussed in \cite{HM2020}.

In this contribution, we investigate some important dynamical properties of poly-PL kinetics systems such as existence and parametrization of complex balanced steady states, capacity for multiple complex balanced steady states in a stoichiometric compatibility class, concentration robustness of some species, and linear stability of complex balanced equilibria.  The primary approach is to specify power law representations of poly-PL kinetic systems in order to apply or mimic existing results on power law systems.  Two approaches are considered in this study:  the \textit{canonical PL-representation} and the \textit{STAR-MSC transformation} \cite{MAGP2019} of the poly-PL system (Section \ref{sec:PYK}).

The main contribution of this paper are the following:
\begin{itemize}
\item \textit{Section 3:} The results of M\"{u}ller and Regensburger \cite{MURE2012,MURE2014} for generalized mass action systems (GMAS) are extended to poly-PL systems. In particular, this work presents partial extensions of the GMAS results to a subset of poly-PL kinetics that are reactant-determined (denoted by PY-RDK systems) and full analogue results for PY-RDK systems with equilibria sets called PL-equilibria and PL-complex balanced equilibria. 

\item \textit{Section 4:} The set of weakly reversible poly-PL systems with zero kinetic reactant deficiency are shown to be PL-complex balanced.

\item \textit{Section 5:} Extensions of the results of Boros et al. \cite{BMR2020} on linear stability of complex balanced equilibria to poly-PL systems are also obtained.

\item \textit{Section 6:} Analogues of the results of Fortun and Mendoza \cite{FM2020} involving \textit{absolute concentration robustness} (ACR) and Lao et al. \cite{LLMM2021} on \textit{balanced concentration robustness} (BCR) for power law kinetic systems to poly-PL kinetic systems are established. Sufficient conditions for concentration robustness in higher deficiency systems are also derived (Section 6.2).


\end{itemize}


\section{Fundamentals of chemical reaction networks and kinetic systems} \label{sec:fund}

We denote by $\mathbb{R}$ and $\mathbb{Z}$ the set of real numbers and integers, respectively. For integers $a$ and $b$, let $\overline{a,b}= \{ j \in \mathbb{Z} | a \leq j \leq b \}.$ We denote the non-negative real numbers by $\mathbb{R}_{\geq0}$, and the positive real numbers by $\mathbb{R}_{>0}$. The sets $\mathbb{R}_{\geq 0}^p$ and $\mathbb{R}_{>0}^p$ are called the \textit{non-negative} and \textit{positive orthants} of $\mathbb{R}^p$, respectively. For $x \in \mathbb{R}^p$, the $i$th coordinate of $x$ is denoted by $x_i$, where $i \in \overline{1,p}$. The \textit{standard basis} for $\mathbb{R}^p$ is the set $\{ \omega_i \in \mathbb{R}^p | i \in\overline{1,p} \}$ and for each $x \in \mathbb{R}^p$ we have the representation $x = \sum_{i=1}^p x_i\omega_i$. Finally, for a vector space $V$, if  $S$ is a non-empty subset of $V$, the span of $S$ is denoted by $\left\langle S \right\rangle$. 

\subsection{Chemical reaction networks basics}
A chemical reaction network (CRN) is a system of interdependent chemical reactions. Each reaction is represented as an ordered pair of vectors, called complexes, of chemical species. 

\begin{definition}
\sloppy A \textbf{chemical reaction network} (CRN) $\mathscr{N}$ is a triple $(\mathscr{S},\mathscr{C},\mathscr{R})$ of three finite sets:
(1) a set $\mathscr{S}= \{X_1, X_2, \dots, X_m \}$ of \textbf{species}, (2) a set $\mathscr{C}$ of \textbf{complexes}, consisting of nonnegative linear combinations of the species, and (3) a set $\mathscr{R} = \{R_1, R_2, \dots, R_r \}\subset \mathscr{C} \times \mathscr{C}$ of \textbf{reactions} such that $(y,y) \notin \mathscr{R}$ for any $y \in \mathscr{C}$, and  for each $y \in \mathscr{C}$, there exists $y' \in \mathscr{C}$ such that either $(y,y') \in \mathscr{R}$ or  $(y',y) \in \mathscr{R}$. We denote the number of species by $m$, the number of complexes with $n$ and the number of reactions by $r$.
\end{definition}

Alternatively, a CRN can be defined as a directed graph whose vertices are embedded as vectors with nonnegative coefficients in a finite dimensional Euclidean space. In this perspective, the vertices correspond to the complexes and the arcs to the reactions. The species can be identified with the standard unit vectors of the Euclidean space. We assume, as is usual in the literature, that i) each species occurs in at least one complex and ii) each complex occurs in at least one reaction.

We use the convention that an element $(y, y') \in \mathscr{R}$ is denoted by $y \rightarrow y' $. In this reaction, we say that $y$ is the \textbf{reactant} complex and $y'$ is the \textbf{product} complex. Connected components of a CRN are called \textbf{linkage classes}, strongly connected components are called \textbf{strong linkage classes}, and strongly connected components without outgoing arcs are called \textbf{terminal strong linkage classes}. We denote the number of linkage classes by $\ell$, that of the strong linkage classes with $s\ell$, and that of terminal strong linkage classes with $t$. A complex is called \textbf{terminal} if it belongs to a terminal strong linkage class. Otherwise, the complex is called \textbf{nonterminal}. A CRN is \textbf{weakly reversible} if $s\ell = \ell$ (i.e., every linkage class is a strong linkage class) and $\bm{t}$\textbf{-minimal} if $t=\ell$ (i.e., every linkage class has one terminal strong linkage class). A CRN is cycle terminal if $n_r=n$.

With each reaction $y\rightarrow y'$, we associate a \textbf{reaction vector} obtained by subtracting the reactant complex $y$ from the product complex $y'$. The \textbf{stoichiometric subspace} $S$ of a CRN is the linear subspace of $\mathbb{R}^m$ defined by
$S := \text{span } \left\lbrace y' - y \in \mathbb{R}^m \mid y\rightarrow y' \in \mathscr{R}\right\rbrace. $
The \textbf{rank} of the CRN  is defined as $s := \dim S$.  Furthermore, if $c \in \mathbb{R}^m_{>0}$, the corresponding \textbf{stoichiometric compatibility class} is defined as the intersection of the coset $c+S$ with $\mathbb{R}^m_{\geq 0}$. The  \textbf{deficiency} of a CRN is the nonnegative integer defined by $\delta = n - \ell - s$. 

\begin{definition}
Let $\mathscr{N}=(\mathscr{S,C,R})$ be a CRN. The \textbf{map of complexes} $\displaystyle{Y: \mathbb{R}^\mathscr{C} \rightarrow \mathbb{R}^\mathscr{S}_{\geq 0}}$ maps the basis vector $\omega_y$ to the complex $ y \in \mathscr{C}$. The \textbf{incidence map} $\displaystyle{I_a : \mathbb{R}^\mathscr{R} \rightarrow \mathbb{R}^\mathscr{C}}$ is the linear map defined by mapping for each reaction $\displaystyle{R_j: y_j \rightarrow y_j' \in \mathscr{R}}$, the basis vector $\omega_j$ to the vector $\omega_{y_j'}-\omega_{y_j} \in \mathscr{C}$. The \textbf{stoichiometric map} $\displaystyle{N: \mathbb{R}^\mathscr{R} \rightarrow \mathbb{R}^\mathscr{S}}$ is defined as $N = Y \circ  I_a$. 
\end{definition}

\subsection{Chemical kinetic system}

By \textit{kinetics} of a CRN, we mean the assignment of a rate function to each reaction in the CRN. It is defined formally as follows \cite{WIUF2013}.
\begin{definition}
A \textbf{kinetics} of a CRN $\mathscr{N}=(\mathscr{S},\mathscr{C},\mathscr{R})$ is an assignment of a rate function $\displaystyle{K_{j}: \Omega_K \to \mathbb{R}_{\geq 0}}$ to each reaction $R_j \in \mathscr{R}$, where $\Omega_K$ is a set such that $\mathbb{R}^{m}_{>0} \subseteq \Omega_K \subseteq \mathbb{R}^{m}_{\geq 0}$. A kinetics for a network $\mathscr{N}$ is denoted by $\displaystyle{K=[K_1,K_2,...,K_r]^\top:\Omega_K \to {\mathbb{R}}^{r}_{\geq 0}}.$ The pair $(\mathscr{N},K)$ is called the \textbf{chemical kinetic system (CKS)}.
\end{definition}

In this paper, in order to study the dynamics of chemical kinetic systems through systems of ordinary differential equations (ODE), we generally assume that the interaction functions are differentiable functions of the variables.

Once a kinetics is associated with a CRN, we can determine the rate at which the concentration of each species evolves at composition $c \in \mathbb{R}^m_{\geq 0}$. 

\begin{definition}\label{def:SFRF}
The \textbf{species formation rate function}  of a chemical kinetic system $(\mathscr{N},K)$ is given by
$$f(c) = NK (c) = \displaystyle\sum_{y_j\rightarrow y'_j \in \mathscr{R}}K_j(c) (y_j'- y_j).$$
\noindent The equation $\text{d}c/\text{d}t=f(c)$ is the \textbf{ODE or dynamical system} of the CKS.  A \textbf{positive equilibrium or steady state} $c^*$ is an element of $\mathbb{R}^m_{>0}$ for which $f(c^*) = 0$.The \textbf{set of positive equilibria} of the system is given by $E_+(\mathscr{N},K)= \{ c^* \in \mathbb{R}^m_0 \mid f(c^*)=0 \}$.
\end{definition}
A chemical kinetic system is \textbf{multistationary} (or has the capacity to admit multiple steady states) if there exist positive rate constants such that the corresponding ODEs system admits at least two distinct stoichiometrically compatible equilibria. Otherwise, it is \textbf{monostationary}.

The \textit{complex formation rate function} is the analogue of the species formation rate function for complexes. 
 
\begin{definition}\label{def:CFRF}
The \textbf{complex formation rate function}  $g: \mathbb{R}^m_{>0} \rightarrow \mathbb{R}^n$ of a chemical kinetic system is the given by
\begin{equation}\label{eq:CFRF}
g(c) = I_a K (c) = \displaystyle\sum_{y_j\rightarrow y'_j \in \mathscr{R}}K_j(c) (\omega_{y_j'}- \omega_{y_j}).
\end{equation} 
\end{definition}

Horn and Jackson \cite{HORNJACK1972} introduced the notion of \textit{complex balancing} in chemical kinetics, which proved to have profound uses in CRNT. Observe from Equation (\ref{eq:CFRF}) that the function $g$ gives the difference between the production and degradation of each complex. Thus, ``complex balancing'' occurs when $g(c)=0$. A chemical kinetic system $(\mathscr{N},K)$ is called \textbf{complex balanced} if it has a complex balanced steady state. The set of positive complex balanced steady states of the system is given by $Z_+(\mathscr{N},K)=\{ c \in \mathbb{R}^m_0 \mid g(c)=0 \}.$ 

Power law kinetics is defined by an  $r \times m$ matrix $F=[F_{ij}]$, called the \textbf{kinetic order matrix}, and vector $k \in \mathbb{R}^r_{>0}$, called the \textbf{rate vector}.  

\begin{definition}
A kinetics $K: \mathbb{R}^m_{>0} \rightarrow \mathbb{R}^r$ is a \textbf{power law kinetics} (PLK) if
$$\displaystyle K_{i}(x)=k_i x^{F_{i,\cdot}} \quad \text{for all } i \in \overline{1,r},$$
with $k_i \in \mathbb{R}_{>0}$ and $F_{ij} \in \mathbb{R}$.  A PLK system has \textbf{reactant-determined kinetics} (of type \textbf{PL-RDK}) if for any two reactions $R_i$, $R_j \in \mathscr{R}$ with identical reactant complexes, the corresponding rows of kinetic orders in $F$ are identical, i.e. $F_{ih}=F_{jh}$ for $h  \in \overline{1,m}$.  Otherwise, a PLK system has \textbf{non-reactant-determined kinetics} (of type \textbf{PL-NDK}).
\end{definition}

An example of PL-RDK is the well-known \textbf{mass action kinetics} (MAK), where the kinetic order matrix is the transpose of the map of complexes $Y$ \cite{FEIN1979}. 

Arceo et al. \cite{AJMSM2015} identified two large sets of kinetic systems, namely the \textbf{complex factorizable} (CF) kinetics and its complement, the \textbf{non-complex factorizable} (NF) kinetics. Complex factorizable kinetics generalize the key structural property of MAK that the species formation rate function decomposes as
$
\text{d}x / \text{d} t= Y \circ A_k \circ \Psi_K,
$
where $Y$ is the map of complexes, $A_k$ is the Laplacian map, and $\Psi_K: \mathbb{R}^\mathscr{S}_{\geq 0} \rightarrow  \mathbb{R}^\mathscr{C}_{\geq 0}$ such that $I_a \circ K(x) = A_k \circ \Psi_K(x)$ for all $x \in \mathbb{R}^\mathscr{S}_{\geq 0}$. In the set of power law kinetics, the complex-factorizable kinetic systems are precisely the PL-RDK systems.

\subsection{Decomposition theory}

We refer to \cite{FML2020} for more details on the concepts and results in decomposition theory.

\begin{definition}
Let $\mathscr{N} =(\mathscr{S}, \mathscr{C}, \mathscr{R})$ be a CRN. A \textbf{covering} of $\mathscr{N}$ is a collection of subsets $\{ \mathscr{R}_1, \mathscr{R}_2,\dots, \mathscr{R}_p \}$ whose union is $\mathscr{R}$. A covering is called a \textbf{decomposition} of $\mathscr{N}$ if the sets $\mathscr{R}_i$ form a partition of $\mathscr{R}$.
\end{definition}

Each $\mathscr{R}_i$ defines a subnetwork $\mathscr{N}_i$ of $\mathscr{N}$ wherein $\mathscr{C}_i$ consists of all complexes occurring in $\mathscr{R}_i$ and $\mathscr{S}_i$ consists of all the species occurring in $\mathscr{C}_i$.

Feinberg  \cite{FEIN1987} identified an important class of network decomposition called {\em independent decomposition}. A decomposition is \textbf{independent} if the stoichiometric subspace $S$ of a network is the direct sum of the subnetworks' stoichiometric subspaces $S_i$ or equivalently, if $s = s_1 + s_2 + \cdots + s_p$. Fortun et al.\cite{FMRL2019} derived a basic property of independent decompositions:

\begin{proposition}\label{prop1}
If $\mathscr{N}=\mathscr{N}_1 \cup \mathscr{N}_2 \cup \cdots  \cup \mathscr{N}_p$ is an independent decomposition, then $
\delta \leq\delta_1 +\delta_2 + \cdots +\delta_p$, where $\delta_i$ represents the deficiency of the subnetwork $\mathscr{N}_i$. 
\end{proposition}

Feinberg \cite{FEIN1987} established the following relationship between the positive equilibria of the ``parent network'' and those of the subnetworks of an independent decomposition:

\begin{theorem}[Feinberg Decomposition Theorem, \cite{FEIN1987}] \label{feinberg theorem}
Let $\{\mathscr{R}_1,\mathscr{R}_2,\dots, \mathscr{R}_p \}$ be a partition of a CRN $\mathscr{N}$ and let $K$ be a kinetics on $\mathscr{N}$. If $\mathscr{N}=\mathscr{N}_1 \cup \mathscr{N}_2 \cup \cdots \cup\mathscr{N}_p$ is the network decomposition generated by the partition  and $E_+(\mathscr{N}_i,K_i)= \{ x \in \mathbb{R}^m_{>0} | N_i K_i(x) = 0 \}$, then \begin{enumerate}
\item[(i)] $ \displaystyle{\bigcap_{i\in \overline{1,p}}} E_+ (\mathscr{N}_i, K_i)  \subseteq E_+ (\mathscr{N}, K)$
\item[(ii)] If the network decomposition is independent, then equality holds.
\end{enumerate}
\end{theorem}

Farinas et al. \cite{FML2020}  introduced the concept of an \textit{incidence independent decomposition}, which naturally complements the independence property. A decomposition of a CRN $\mathscr{N}$ is \textbf{incidence independent} if  the image of the incidence map $I_a$ of $\mathscr{N}$ is the direct sum of the images of the incidence maps of the subnetworks. It follows that the dimension of the image of the incidence map $I_a$ equals the sum of the dimensions of the subnetworks' incidence maps. That is, $n-\ell = \sum (n_i - \ell_i)$. The linkage classes form the primary example of an incidence independent decomposition, since $n = \sum n_i$ and $\ell = \sum \ell_i$. 

The following result is the analogue of Proposition \ref{prop1} for incidence independent decomposition.

\begin{proposition}[Prop. 7, \cite{FML2020}]
\label{prop:incidenceindep}
Let $\mathscr{N}=\mathscr{N}_1 \cup \mathscr{N}_2 \cup \cdots \cup \mathscr{N}_p$ be an incidence independent decomposition. Then $\delta \geq \delta_1 +\delta_2 + \cdots + \delta_p$.
\end{proposition}

$\mathscr{C}$-decompositions form an important class of incidence independent decompositions:

\begin{definition}\label{def:Cdecomp}
A decomposition $\mathscr{N} = \mathscr{N}_1 \cup \mathscr{N}_2 \cup \cdots \cup \mathscr{N}_p$ with $\mathscr{N}_i =(\mathscr{S}_i, \mathscr{C}_i, \mathscr{R}_i)$ is a \textbf{$\mathscr{C}$-decomposition} if $\mathscr{C}_i \cap \mathscr{C}_j = \emptyset$ for $i \neq j$.
\end{definition}

A $\mathscr{C}$-decomposition partitions not only the set of reactions but also the set of complexes. The primary examples of  $\mathscr{C}$-decomposition are the linkage classes.

The following result shows the relationship between the set of incidence independent decompositions and the set of complex balanced equilibria of any kinetic system.  It is the precise analogue of Feinberg's result (Theorem \ref{feinberg theorem}). 

\begin{theorem}[Theorem 4, \cite{FML2020}]
\label{th:Z}
Let $\mathscr{N}=(\mathscr{S}, \mathscr{C}, \mathscr{R})$ be a CRN and $\mathscr{N}_i =(\mathscr{S}_i, \mathscr{C}_i, \mathscr{R}_i)$ for $i\in \overline{1,p}$ be the subnetworks of a decomposition. Let $K$ be any kinetics, and $Z_+ (\mathscr{N},K)$ and $Z_+ (\mathscr{N}_i, K_i)$  be the set of  complex balanced equilibria of $\mathscr{N}$ and $\mathscr{N}_i$, respectively. Then
\begin{enumerate} 
\item[(i)] $\displaystyle{\bigcap_{i\in \overline{1,p}}} Z_+ (\mathscr{N}_i, K_i) \subseteq Z_+ (\mathscr{N}, K)$
\end{enumerate}
 If the decomposition is incidence independent, then 
 \begin{enumerate}
\item[(ii)] $Z_+ (\mathscr{N}, K)= \displaystyle{\bigcap_{i\in \overline{1,p}}} Z_+ (\mathscr{N}_i, K_i)$
\item[(iii)] $ Z_+ (\mathscr{N}, K) \neq \emptyset$ implies $Z_+ (\mathscr{N}_i, K_i) \neq \emptyset$ for each $i\in \overline{1,p}$.
\end{enumerate}
\end{theorem}

The converse statement of Theorem \ref{th:Z} (iii) holds for a subset of incidence independent decompositions with any kinetics:

\begin{proposition}[Theorem 5, \cite{FML2020}]
Let $\mathscr{N}=\mathscr{N}_1 \cup \mathscr{N}_2 \cup \cdots \cup \mathscr{N}_p$ be a weakly reversible $\mathscr{C}$-decomposition of a chemical kinetic system $(\mathscr{N},K)$. If $Z_+(\mathscr{N}_i, K_i) \neq \emptyset$ for each $i \in \overline{1,p}$, then $Z_+(\mathscr{N},K)\neq \emptyset$. 
\end{proposition}
\section{Complex balanced poly-PL kinetic systems}\label{sec:PYK}

Poly-PL systems were introduced by Talabis et al. in \cite{TMMNJ2020} and applied to study evolutionary games with replicator dynamics and polynomial payoff functions. A computational approach to determine their multistationarity was developed by Magpantay et al. \cite{MAGP2019}. In this section, we first collect concepts, techniques and results from these papers and then apply them to the extension of the results of M\"{u}ller and Regensburger on Generalized Mass Action Systems (GMAS) to classes of poly-PL systems.

\subsection{Poly-PL systems: canonical PL-representation and STAR-MSC transformation} \label{sec:polyPL}

We recall the definition and introduce some additional notation for a PY-RDK system from Talabis et al. \cite{TMMNJ2020}:

\begin{definition}
A kinetics $K:\mathbb{R}^m_{>0} \rightarrow \mathbb{R}^r$ is a \textbf{poly-PL kinetics} if 
\begin{equation}
K_i(x) = k_i \left( a_{i,1}x^{F_{i,1}}+a_{i,2}x^{F_{i,2}}+\cdots + a_{i,{h_i}}x^{F_{i,{h_i}}} \right) \quad \text{for } i \in \overline{1,r}
\end{equation}
written in lexicographic order with $k_i \in \mathbb{R}_{>0}$, $a_{i,j} \in \mathbb{R}_{\geq 0}$, $F_{i,j} \in \mathbb{R}^m$, and $j \in \overline{1,h_i}$ (where $h_i$ is the number of terms in reaction $i$).  Power-law kinetics is defined by the $r \times m$ \textbf{kinetic order matrices} $\left[ F_{i,j} \right]$, the \textbf{rate vector} $ \left[ k_i \right] \in \mathbb{R}^r_{>0}$, and the \textbf{poly-rate vectors} $\left[ a_{i,j} \right]\in \mathbb{R}^r_{>0} $. If $h=\max h_i$, we normalize the length of each kinetics to $h$ by replacing the last term with  $\left( h- h_i + 1\right)$ terms with $ \frac{1}{h-h_i +1} x^{F_{i,h_i}}$. We call this the \textbf{canonical representation} of a poly-PL kinetics. 
\end{definition}

\noindent Henceforth, for each $j \in \overline{1,h}$, we set $\mathcal K_j (x) := k_i a_{i,j} x^{F_{i,j}}$ where $i \in \overline{1,r}$. Furthermore, we use the notation $\{ \mathcal K_j \}$ (where $j \in \overline{1,h}$) or $K =\mathcal K_1 + \mathcal K_2 +\cdots + \mathcal K_h $ to denote the PL-representation of a poly-PL kinetics. \\

We extend the definition of PL-RDK kinetics to poly-PL kinetics. 

\begin{definition}
A poly-PL kinetic system $(\mathscr{N},K)$ is said to be have \textbf{reactant-determined kinetics} if and only if for each $(\mathscr{N},\mathcal K_j)$ obtained from the PL-representation $\{\mathcal K_j \}$ of the system, $(\mathscr{N},\mathcal K_j)$ is PL-RDK and $a_{i,j} = a_{i',j}$ for any two reactions $R_i,R_{i'}$ of a branching node. We denote such poly-PL system by \textbf{PY-RDK}. 
\end{definition}

The \textbf{S-invariant termwise addition of reactions via maximal stoichiometric coefficients (STAR-MSC)} method is based on the idea to use the maximal stoichiometric coefficient (MSC) among the complexes in the CRN to construct reactions whose reactant complexes and product complexes are different from existing ones. This is done by uniform translation of the reactants and products to create a ``replica'' of the CRN. The method creates $h - 1$  replicas of the original network  and hence its transform, $\mathscr{N}^*$ becomes the union (in the sense of \cite{GHMS2020}) of the replicas and the original CRN.

\sloppy We now describe the STAR-MSC transformation. Since the domain of definition of a poly-PL kinetics is $\mathbb{R}^m_{>0}$, all $x=(X_1,X_2,\dots, X_m)$ are positive vectors. Let $\displaystyle{M = 1 + \max \{ y_i  \mid y \in \mathscr{C} \}}$, where the second summand is the maximal stoichiometric coefficient. 

For any positive integer $z$, define the vector $\bm{\underline{z}}$ to be the vector $(z,z,\dots, z) \in \mathbb{R}^m$. For each complex $y \in \mathscr{C}$, form the $(h-1)$ complexes
$$
y + \bm{\underline{M}}, y + \bm{\underline{2M}}, \dots, y + \bm{\underline{(h-1)M}}.
$$
Each of these complexes are different from each other and from all existing complexes. Further properties of the STAR-MSC transformation are discussed in \cite{MAGP2019}.

\begin{example}\label{example1}
Consider the following weakly reversible network with two species $X$ and $Y$.  The kinetics are also given below.
\begin{equation}
\nonumber
\begin{tikzpicture}[baseline=(current  bounding  box.center)]
\tikzset{vertex/.style = {draw=none,fill=none}}
\tikzset{edge/.style = {bend left,->,> = latex', line width=0.20mm}}
\node[vertex] (1) at  (0,0) {$X$};
\node[vertex] (2) at  (4,0) {$2X+Y$};
\node[vertex] (3) at  (0,-2.5) {$Y$};
\node[vertex] (4) at  (4,-2.5) {$X+2Y$};
\draw [edge]  (1) to["$k_1$"] (2);
\draw [edge]  (2) to["$k_2$"] (1);
\draw [edge]  (3) to["$k_3$"] (4);
\draw [edge]  (4) to["$k_4$"] (3);
\end{tikzpicture}
\quad
K(\textbf{X})=\left[ 
\begin{array}{ccc}
	k_1 x \\
	k_2 y \\
	k_3 \frac{x}{1+x} \\
	k_4 \frac{y}{1+y}  \\	
\end{array}
 \right]
\end{equation}

\noindent $K_3$ and $K_4$ are of Hill type kinetics. Using the transformation of Hernandez and Mendoza \cite{HM2020}, the associated kinetics will be a poly-PL kinetics given by:

\begin{equation}
\nonumber
K_{PY}(\textbf{X})=\left[ 
\begin{array}{ccc}
	k_1 x(1+x)(1+y) \\
	k_2 y(1+x)(1+y) \\
	k_3 x(1+y) \\
	k_4 y(1+x)  \\	
\end{array}
 \right].
\end{equation}
\end{example}

\noindent The canonical PL-representation of a poly-PL kinetics is 
\begin{equation}
\nonumber
K_{PY}(\textbf{X})=\left[ 
\begin{array}{c}
	k_1 (x+x^2+xy+x^2y) \\
	k_2 (y+xy+y^2+xy^2) \\
	k_3 (xy+\frac{1}{3}x+\frac{1}{3}x+\frac{1}{3}x) \\
	k_4 (xy+\frac{1}{3}y+\frac{1}{3}y+\frac{1}{3}y)  \\	
\end{array}
 \right].
\end{equation}

\begin{example}\label{example2}
Consider the following Michaelis-Menten kinetic system from Enzyme biology with set of biochemical species $S = \left\{S_1, S_2, S_3, S_4\right\}$, set of complexes $\left\{S_1+S_2, S_1+S_3, S_4\right\}$, and reactions:

\begin{equation}
\nonumber
\begin{tikzpicture}[baseline=(current  bounding  box.center)]
\tikzset{vertex/.style = {draw=none,fill=none}}
\tikzset{edge/.style = {bend left,->,> = latex', line width=0.20mm}}
\node[vertex] (1) at  (0,0) {$S_1+S_2$};
\node[vertex] (2) at  (2,0) {$S_4$};
\node[vertex] (3) at  (4,0) {$S_1+S_3$};
\draw [edge]  (1) to["$k_1$"] (2);
\draw [edge]  (2) to["$k_2$"] (1);
\draw [edge]  (2) to["$k_3$"] (3);
\draw [edge]  (3) to["$k_4$"] (2);
\end{tikzpicture}
\end{equation}

\noindent This network assumes the Michaelis–Menten enzyme mechanism, in which a substrate $S_2$ is modified into a substrate $S_3$ through the formation of an intermediate $S_4$. The reaction is catalyzed by an enzyme $S_1$. The kinetic systems associated to the ODE system are:
$$K_{MAK}(\textbf{X}) = \left[ 
\begin{array}{c}
k_1 S_1 S_2\\
k_2  S_3 \\
k_3  S_3 = v_f \\
k_4 S_1 S_3 = v_b \\
\end{array}
 \right] \text{ and } K_{PQK}(\textbf{X}) = \left[ 
\begin{array}{c}
k_1 S_1 S_2\\
k_2  S_3 \\
k_1  \frac{S_2}{1+k_3 S_2 + S_4} \\
k_2  \frac{S_4}{1+k_3 S_2 + S_4} \\
\end{array}
 \right].$$

The kinetics $K_{PQK}(\textbf{X})$ is taken from the collection of enzymatic reactions in Appendix 2 of \cite{GABOR2015}. The forward rate $v_f$ and the backward rate $v_b$ are replaced by rational expressions $k_1  \frac{S_2}{1+k_3 S_2 + S_4}$ and $k_2  \frac{S_4}{1+k_3 S_2 + S_4}$, respectively. \\

Using the transformation in Section 6.1 of \cite{HM2020}, the associated kinetics will be a poly-PL kinetics in the form:

\begin{equation}
\nonumber
K_{PY}(\textbf{X})=\left[  
\begin{array}{c}
k_1 (S_1 S_2+k_3 S_1 S_2^2 + S_1 S_2 S_4) \\
k_2 (S_3+k_3 S_2 S_3 + S_3 S_4) \\
k_1  S_2 \\
k_2  S_4 \\
\end{array}
 \right].
\end{equation}

The canonical PL-representation of a poly-PL kinetics is 
\begin{equation}
\nonumber
K_{PY}(\textbf{X})=\left[ 
\begin{array}{c}
k_1 (S_1 S_2+k_3 S_1 S_2^2 + S_1 S_2 S_4) \\
k_2 (S_3+k_3 S_2 S_3 + S_3 S_4) \\
k_1  (\frac{1}{3}S_2 + \frac{1}{3}S_2 + \frac{1}{3}S_2) \\
k_2  (\frac{1}{3}S_4 + \frac{1}{3}S_4 + \frac{1}{3}S_4) \\
\end{array}
 \right].
\end{equation}

\end{example}


\subsection{Extension of GMAS theory to PY-RDK systems}

A kinetics $K$ on a weakly reversible network $\mathscr{N}$ has \textbf{conditional complex balancing (CCB)} if there exist rate constants such that $(\mathscr{N},K)$ is complexed balanced, i.e., $Z_+(\mathscr{N} ,K) \neq \emptyset$. A kinetics is \textbf{unconditionally complex balanced (UCB)} on a weakly reversible network if for all rate constants, $Z_+(\mathscr{N} ,K) \neq \emptyset$.

Poly-PL kinetic systems belong to a large set of chemical kinetics called \textbf{rate constant-interaction map decomposable (RID) kinetics}. RID kinetic systems are chemical systems with constant rates \cite{NEML2019}. In Appendix \ref{append:A}, we define this system formally and we derive the Conditional Complex Balancing (CCB) property for any poly-PL kinetics on a weakly reversible network from a general result on RID kinetics. 

We recall, from \cite{MURE2014}, that a \textbf{generalized mass action system (GMAS)} is a 4-tuple $(G, y, \widetilde{y}, k)$ where $G$ is a digraph with $n$ vertices and $\ell$ connected components, $y$ maps the complexes into $\mathbb{R}^m$, $\widetilde{y}$ maps the reactant complexes to $\mathbb{R}^m$ and $k$ is a vector of rate constants. The image of $y$ can be viewed as the (stoichiometric) complexes of a CRN. The image of $\widetilde{y}$ is the set of \textbf{kinetic complexes}. It was shown in \cite{TAM2018} that a PL-RDK system corresponds to a GMAS with an injective map $y$. 

The cornerstone of the theory of M\"{u}ller-Regensburger for GMAS systems \cite{MURE2012} is the concept of kinetic deficiency $\widetilde{\delta}$, which in turn is based on the concept of a kinetic order subspace $\widetilde{S}$. The kinetic deficiency is defined as $\widetilde{\delta}=n-l-\dim ( \widetilde{S})$.
Recall that for a weakly reversible GMAS system, the kinetic order subspace $\widetilde{S}$ is the span $\langle \widetilde{y}'-\widetilde{y}\rangle$ where $\widetilde{y}',\widetilde{y}$ are kinetic complexes of $y\rightarrow y'$ in $\mathbb{R}^\mathscr{S}$. 

We extend the concept of kinetic order subspace to a PY-RDK system as follows: to each reaction, we assign the $h$-vector of kinetic complexes of $(\mathscr{N}, \mathcal K_j)$ in the direct sum of $\mathbb{R}^m$ $h$-times. Since the system is complex factorizable, this is also map from the set of reactant complexes of $\mathscr{N}$.  We define $\widetilde{S}$ as follows:

\begin{definition}
The \textbf{kinetic order subspace} of a poly-PL system is the direct sum $\widetilde{S}=\widetilde{S}_1 + \cdots + \widetilde{S}_h$, considered as a subspace of the direct sum of $\mathbb{R}^m$ $h$-times. Here, $\widetilde{S}_j$ is the kinetic order subspace of the GMAS system $(\mathscr{N}, \mathcal K_j)$ (where $j \in \overline{1,h}$). The \textbf{dimension} of the kinetic order subspace, $\widetilde{s}_j$, and the \textbf{kinetic deficiency}, $\widetilde{\delta}$, of the poly-PL system $(\mathscr{N},K)$ are defined as $\widetilde{s}=\widetilde{s}_1 + \cdots + \widetilde{s}_h$ and
$\widetilde{\delta}= \widetilde{\delta}_1 + \cdots + \widetilde{\delta}_h$, respectively. Here, $\widetilde{s}_j = \dim \widetilde{S}_j$ and 
$\widetilde{\delta}_j =  n - l - \widetilde{s}_j$ (for each $j \in \overline{1,h}$).
\end{definition}

We present partial extensions of the GMAS results to PY-RDK systems. First, we show that zero kinetic deficiency is sufficient for unconditional complex balancing:

\begin{proposition}
Let $(\mathscr{N}, K)$ be a weakly reversible PY-RDK system with $\widetilde{\delta}=0$. Then $Z_+(\mathscr{N},K) \neq \emptyset$ for any set of rate constants.
\end{proposition}

\textit{Proof.}
\sloppy Let $(\mathscr{N}^*,K^*)$ be the STAR-MSC transform of $(\mathscr{N},K)$. Since $\mathscr{N}$ is weakly reversible, each $\mathscr{N}^*_j$ is weakly reversible. Moreover, $\widetilde{\delta}=0$ implies that  $\widetilde{\delta}_j=0$. By the GMAS unconditional complex balancing, we have $Z_+(\mathscr{N},\mathcal K_j) \neq \emptyset$ and hence, $Z_+(\mathscr{N}_j^*,K_j^*) \neq \emptyset$. Since these sets form a $\mathscr{C}$-decomposition of $(\mathscr{N}_j^*,K_j^*)$, we have $\emptyset \neq Z_+(\mathscr{N}_j^*,K_j^*)=\bigcap Z_{+}(\mathscr{N},\mathcal K_j)$. Therefore, $Z_+(\mathscr{N},K) \neq \emptyset$ for any set of rate constants. $\blacksquare$

The poly-PL systems in Examples \ref{example1} and \ref{example2} are weakly reversible PY-RDK with $\tilde{\delta}=0$. Hence, these systems are unconditionally complexed balanced. The following corollary illustrates the usefulness of the above extension:

\begin{corollary}\label{cor:gmak1}
Let $(\mathscr{N}, K)$ be a weakly reversible PY-RDK system. If its canonical PL-representation $\{ \mathcal K_j \}$ consists of PL-RDK systems with zero kinetic deficiency, then $Z_+(\mathscr{N}, K) \neq \emptyset$ for any set of rate constants.
\end{corollary}


The STAR-MSC transform $(\mathscr{N}^*, K^*)$ as a PL-RDK system also has a kinetic deficiency, which we denote by $\widetilde{\delta}^*$. We document an interesting relation to $\widetilde{\delta}$ as well as two other bounds in the following proposition:

\begin{proposition}
Let $(\mathscr{N}^*, K^*)$ be the STAR-MSC transform of $(\mathscr{N}, K)$. Then
\begin{enumerate}
\item [i.] $\widetilde{\delta}^* \leq h(n-\ell)-1$.
\item [ii.] If $h \geq m$, then $\widetilde{\delta}^* \geq \widetilde{\delta}$.
\item [iii.] If $\mathscr{N}$ is open (i.e., non-mass conserving), then $\widetilde{\delta}^* \geq \delta^*$.
\end{enumerate}
\end{proposition}

\textit{Proof.}
For the kinetic rank $\widetilde{s^*}$ we have $ 1 \leq \widetilde{s^*} \leq m$. Hence, $h(n-\ell)-m \leq \widetilde{\delta^*} \leq h(n-\ell)-1$, which shows (i). For (ii), note that $h(n-\ell)-m= h(n-\ell)- \widetilde{s} + (\widetilde{s} -m)$ so that $\widetilde{\delta^*}  - \widetilde{\delta} \geq \widetilde{s} - m \geq h - m$, since $\widetilde{s}_j \geq 1$. The assumption that $h \geq m$ shows the claim. For (iii), by assumption, $s=m$ and $s=s^*$ so that the claim follows as well.
$\blacksquare$

We introduce equilibria subsets to which the full analogues of the parametrization and monostationarity results for GMAS systems can be extended. 

\begin{definition}
The set of \textbf{PL-equilibria} $E_{+,\text{PL}} (\mathscr{N},K)$ and of \textbf{PL-complex balanced equilibria}  $Z_{+,\text{PL}} (\mathscr{N},K)$ of a poly-PL kinetic system are defined  as 
$
E_{+,\text{PL}} (\mathscr{N},K) :=\bigcap E_+(\mathscr{N},\mathcal K_j)$ and $Z_{+,\text{PL}} (\mathscr{N},K) :=\bigcap Z_+(\mathscr{N}, \mathcal K_j).
$
\end{definition}



A PYK system is called PL-complex balanced if $\emptyset \neq Z_+(\mathscr{N},K)=Z_{+,\text{PL}}(\mathscr{N},K)$. The next result is the parametrization of a subset of the set of complex balanced equilibria:

\begin{proposition}\label{prop:gmak2}
For any PY-RDK system with complex balanced equilibrium $c^* \in
Z_{+,\text{PL}} (\mathscr{N},K)$, 
\begin{equation}
Z_{+,\text{PL}} (\mathscr{N},K) = \{ c \in \mathbb{R}^m_{>0} \mid \ln (c) - \ln (c^*) \in (\widetilde{S})^\perp \}.
\end{equation}
\end{proposition}

\textit{Proof.}
By assumption, for each $j \in \overline{1,h}$, $c^* \in Z_{+} (\mathscr{N},\mathcal K_j)$. By GMAS parametrization (Prop. 2.21 of \cite{MURE2012}), we have $Z_{+} (\mathscr{N},\mathcal K_j) = \{ c \in \mathbb{R}^m_{>0} \mid \ln (c) - \ln (c^*) \in (\widetilde{S}_j)^\perp \}$. This is equivalent to
\begin{align*}
 Z_{+,\text{PL}} (\mathscr{N},K) & = \bigcap \left\lbrace c \in \mathbb{R}^m_{>0} \mid \ln (c) - \ln (c^*) \in (\widetilde{S}_j)^\perp \right\rbrace \\
 &= \left\lbrace c \in \mathbb{R}^m_{>0} \mid \ln (c) - \ln (c^*) \in \bigcap (\widetilde{S}_j)^\perp \right\rbrace \\
 &= \left\lbrace c \in \mathbb{R}^m_{>0} \mid \ln (c) - \ln (c^*) \in (\widetilde{S})^\perp \right\rbrace \text{ since } \bigcap (\widetilde{S}_j)^\perp = \left(\sum \widetilde{S}_j \right)^\perp. 
\end{align*} 
$\blacksquare$

Recall that the \textbf{sign function} of a real vector assigns to each coordinate either $-1$, $1$ or $0$ depending on whether the coordinate is negative, positive or zero.  When applied to a kinetic vector in
$\mathbb{R}^m \times \cdots \times \mathbb{R}^m$, it can be viewed both as a ``long'' vector $(\text{sign}_1, \dots, \text{sign}_h)$ or a vector of $m$-tuples, since the sign function does not depend on the ambient vector space. Finally, we obtain the criterion for monostationarity  in a restricted form:

\begin{proposition}\label{prop:gmak3}
A weakly reversible PY-RDK  system has at most one complex balanced PL-equilibrium in a stoichiometric class if and only if $\text{\em sign} (S) \cap \text{\em sign} (\widetilde{S})^\perp = 0$.
\end{proposition}

\textit{Proof.}
We have 
$$
\text{sign} (\widetilde{S})^\perp = (\text{sign} (\widetilde{S}_1)^\perp,\dots, \text{sign} (\widetilde{S}_h)^\perp) \text{ and } 
\text{sign} (S)  =(\text{sign}(S), \dots, \text{sign} (S)).
$$
For the forward implication, let $c^{*}$ be a PL-equilibrium. It is contained in each $Z_{+} (\mathscr{N},\mathcal K_j)$ where for each $j \in \overline{1,h}$, the sign conditions hold. This results in the $hm$ 0-vector for the overall sign condition. For the backward direction, if the intersection is the $hm$ 0-vector, we have $h$ $m$-tupled zero vectors for the intersections of $\text{sign} (S) \cap \text{ sign} (\widetilde{S_j})^\perp$. If we consider $S$ and $S_j$'s as subspaces of $\text{Im } Y^*$ (i.e. the map of complexes of $\mathscr{N}^*$ in a single $\mathbb{R}^m$), we can conclude, since the sign values do not change, that $(\mathscr{N}^*,K^*)$ is a weakly reversible PL-RDK system with at most one equilibria in every stoichiometric class. Since $Z_+(\mathscr{N}^*,K^*)=Z_{+,\text{PL}}(\mathscr{N},K)$ and $S^* = S$, we obtain the claim. 
$\blacksquare$

\begin{proposition}\label{prop:mure}
\sloppy If a weakly reversible PY-RDK system is PL-complex balanced, i.e., has the
property $Z_+(\mathscr{N},K)=Z_{+,\text{PL}}(\mathscr{N},K)$, then we have:
\begin{enumerate}
\item[i.] \textbf{\em Unconditional complex balancing}: $\widetilde{\delta}=0$ if and only if $Z_+(\mathscr{N},K) \neq \emptyset$ for any set of rate constants
\item[ii.] \textbf{\em Parametrization}: $Z_{+} (\mathscr{N},K) = \{ c \in \mathbb{R}^m_{>0} \mid \ln (c) - \ln (c^*) \in (\widetilde{S})^\perp \}$
\item[iii.]  \textbf{\em Monostationarity criterion}: The poly-PL kinetic system $(\mathscr{N},K)$ is monostationary if and only if $\text{\em sign} (S) \cap \text{\em sign} (\widetilde{S})^\perp = 0$.
\end{enumerate}
\end{proposition}

\section{Weakly reversible poly-PL systems with zero kinetic reactant deficiency}
\label{sec:weakrevpoly}

In this Section, we present the set PY-TIK of poly-PL systems with zero kinetic reactant deficiency as an interesting example of PL-complex balanced systems with zero kinetic deficiency. We first review briefly the results in the special case $h = 1$, i.e. power law systems. We then present the Zero Kinetic Reactant Deficiency (ZKRD) Theorem for PY-TIK, from which the PL-complex balancing property of PY-TIK systems follows. We conclude with a discussion of particular aspects of PY-TIK systems arising from the extension of GMAS properties to PL complex balanced systems presented in Section 3.

\subsection{A brief review of concepts and results for $h = 1$ (PL-TIK systems)}

When studying PL-RDK systems, it is more convenient to use the $m\times n_r$ T-matrix rather than the $r \times m$ kinetic order matrix: $T_{i, \rho(q)} := F^T_{q,i},$ where $i = 1,\dots,m, q = 1,\dots,r$ and $\rho$ is the reactant map. PL-TIK systems were first identified by Talabis et al \cite{TAM2018} as a subset of PL-RDK systems whose T matrix columns in each linkage class were linearly independent (PL-LLK, s. Figure 1). Their defining property was the maximal rank of their augmented T-matrix, e.g. the $(m\times l) \times n_r$ matrix $\widehat{T}$ formed by adding the characteristic function of reactants per linkage class to the T-matrix. 

\begin{figure}[h!]
    \centering
    \includegraphics[width=10cm]{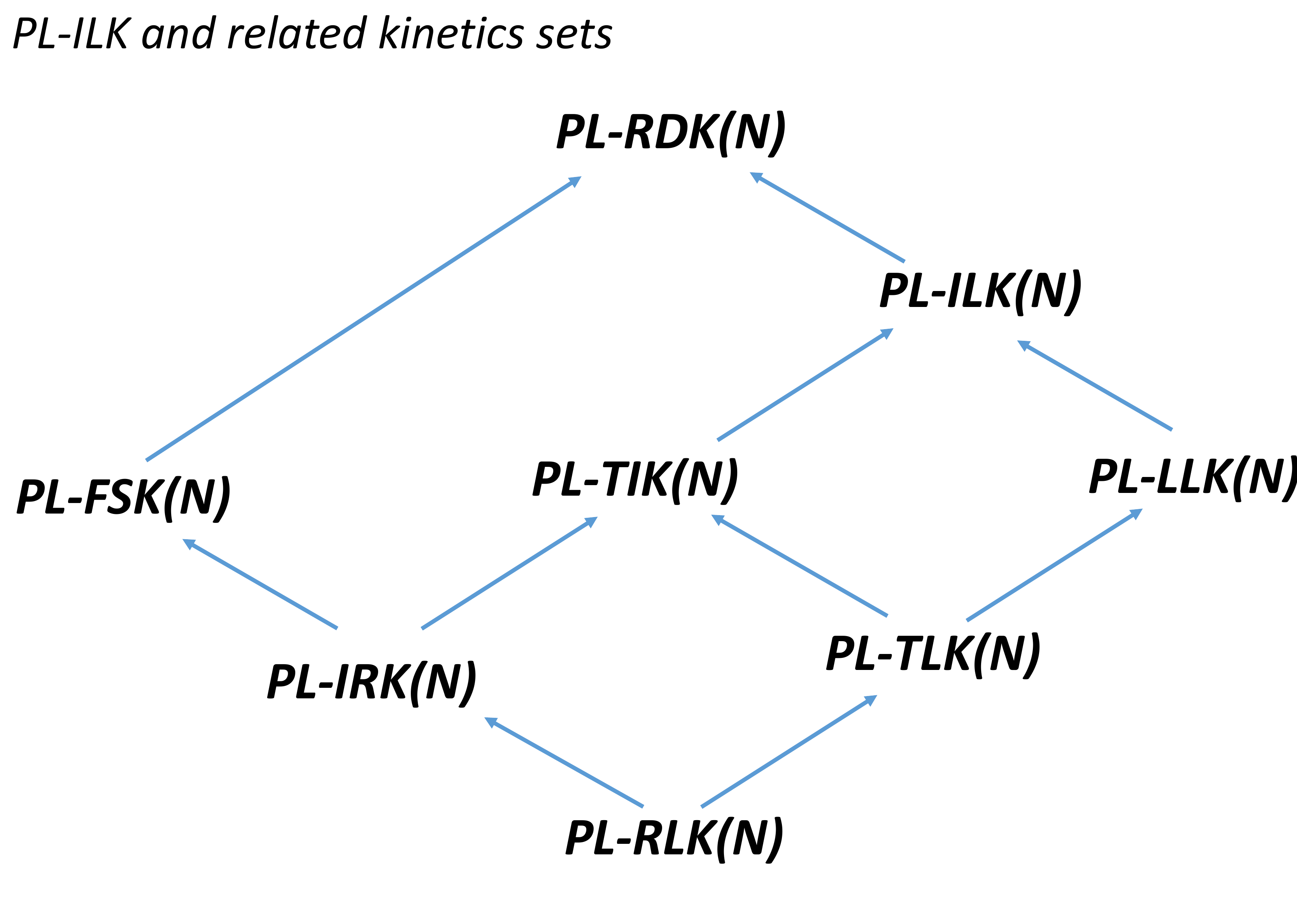}
    \caption{PL-ILK and related kinetic sets}
    \label{fig:figPLTIK}
\end{figure}

Analogues of the Deficiency Zero Theorem (DZT) and the Deficiency One Theorem (DOT) for mass action systems were derived. While the DZT is a special case of a general result by M\"{u}ller and Regensburger for GMAS \cite{MURE2012}, the DOT extension, to our knowledge, remains unique to date. In parallel, Arceo et al \cite{AJLM2017} studied the relations of reactant subspaces of CRN to kinetics and by introducing the concept of reactant deficiency showed similarities to work on PL-TIK. After a special case was derived as a “Weak Reversibility Theorem” by Mendoza et al \cite{MTJ2018}, the full Zero Kinetic Reactant Deficiency (ZKRD) Theorem was established by Talabis et al \cite{TMJ2019} where the new characterization of PL-TIK via the concept of kinetic reactant deficiency was introduced.

\subsection{The Zero Kinetic Reactant Deficiency (ZKRD) Theorem for PY-TIK implies PL-complex balancing}

Let $(\mathscr{N}, K)$ be a PY-RDK system with canonical PL-representation $K = \mathcal K_1 +\dots + \mathcal K_h.$ We now formally define the kinetic reactant deficiency of $(\mathscr{N}, K)$.

For each kinetic order matrix $F_j$ ($\forall j \in \overline{1,h}$), we define the $m \times n$ matrix $\widetilde{Y}_j$ defined as:

\[
    (\widetilde{Y}_j)_{ij}= 
\begin{cases}
    (F_{j})_{ri},& \text{if}\; j\; \text{is a reactant complex of reaction}\; r \\
		0,              & \text{otherwise}
\end{cases}
\]

\begin{definition}
The $m \times n_r$ \textbf{poly T-matrix} $T_j$ ($\forall j \in \overline{1,h}$) is the truncated $\widetilde{Y}_j$ where the non-reactant columns are deleted. Define the $n_r \times l$ matrix $L = \left[ e_1,e_2,...,e_l \right]$ where $e^i$ is a characteristic vector for linkage class $\mathscr{L}^i$. The block matrix $\widehat{T}_j \in \mathbb{R}^{(m+l)\times n_{r}}$ ($\forall j \in \overline{1,h}$) is defined as 
\begin{equation}
\widehat{T}_j=\left[ 
\begin{array}[center]{c} T_j \\
L^{\top} \\
\end{array} \right].
\end{equation}
\end{definition}

\begin{definition}
The block matrix $\widehat{T} \in \mathbb{R}^{h \cdot (m+l)\times h \cdot n_{r}}$ is defined as
\begin{equation}
\widehat{T}=\left[ 
\begin{array}{ccc}
	\widehat{T}_1	& & 0 \\
		&\ddots & \\
		0 & & \widehat{T}_h \\
\end{array}
 \right]
\end{equation}
\end{definition}

\begin{definition}
\label{kinetic reactant deficiency}
Let $\mathscr{N}$ be a network with $n_r$ reactant complexes and $K$ a poly-PL kinetics with poly T-matrices $T_1$, ..., $T_h$. If $\widehat{q}=\text{rank} (\widehat{T})$, then the \textbf{kinetic reactant deficiency} $\widehat{\delta}$ is defined as 
\begin{equation}
\widehat{\delta}= h \cdot n_r - \widehat{q}.
\end{equation}
\end{definition}

\begin{definition}
The set PY-TIK consists of all complex factorizable poly-PL kinetics with zero kinetic reactant deficiency.
\end{definition}

We have the following basic relationships between the kinetic reactant properties of the PY-TIK system and its PL-representation systems. Note that the block matrices $\widehat T_j$ is that of the PLK system $(\mathscr N, \mathcal K_j).$ Set $\widehat q_j$ and $\widehat\delta_j$ be its kinetic reactant rank and kinetic reactant deficiency, respectively.

\begin{proposition}
Let $(\mathscr N, K)$ be a poly-PL system with canonical PL-representation $K = \mathcal K_1 +\dots + \mathcal K_h.$ Then
\begin{enumerate}
    \item [(i)]	$\widehat\delta=\widehat\delta_1+\dots+\widehat\delta_h$; 
    \item [(ii)] $(\mathscr N, K)$ is PY-TIK if and only if each $(\mathscr N, \mathcal K_j)$ is PL-TIK.
\end{enumerate}
\end{proposition}

\textit{Proof.}
(i) follows directly from the definition of the block matrix $\widehat T$ as direct sum of the $\widehat T_j$, which implies $\widehat\delta=\widehat\delta_1+\dots+\widehat\delta_h.$ (ii) follows from (i) directly, too.
$\blacksquare$

We recall from \cite{TMMNJ2020} the complex equilibria existence parametrization and uniqueness statements of the Zero Kinetic Reactant Deficiency (ZRKD) Theorem for PY-TIK:

\begin{theorem}
Let $(\mathscr{N},K)$ be a PY-TIK systems, that is, $(\mathscr{N},K)$ has a complex factorizable kinetics and $\widehat{\delta}=0$. Then $\mathscr{N}$ is weakly reversible if and only if $Z_+(\mathscr{N}, K) \neq \emptyset$.
\end{theorem}

\begin{theorem}
Let $(\mathscr{N},K)$ a weakly reversible poly-PL kinetic system with poly T-matrices $T_1$, ..., $T_h$ and $\widehat{\delta}=0$. Consider an arbitrary poly T-matrix $T_k$.
\label{thm:para}
\begin{enumerate}
\item [(i)] if $Z_+(\mathscr{N}, K) \neq \emptyset$ and $x^* \in Z_+(\mathscr{N}, K)$ then $$Z_+(\mathscr{N},K)=\left \{x\in \mathbb{R}^{m}_\geq \middle| \log(x)-\log(x^*) \in (\widetilde{S}_j)^{\perp} \right \}.$$
\item [(ii)] if $Z_+(\mathscr{N}, K) \neq \emptyset$ then $|Z_+(\mathscr{N},K) \cap Q_j|=1$ for each positive kinetic reactant flux class $Q_j$.
\end{enumerate}
\end{theorem}

\begin{corollary}
A weakly reversible PY-TIK system is PL-complex balanced.
\end{corollary}

\textit{Proof.}
Since $(\mathscr N, \mathcal K_j$) is a weakly reversible PL-TIK, according to the ZKRD Theorem for PL-TIK \cite{MTJ2018}, it is complex balanced and $Z_+(\mathscr{N},\mathcal K_j)=\left \{x\in \mathbb{R}^{\mathscr S}_\geq \middle| \log(x)-\log(x^*) \in (\widetilde{S}_j)^{\perp} \right \}.$ The intersection $\cap Z_+(\mathscr N, \mathcal K_j)$ is, according to Theorem \ref{thm:para} the intersection of $Z_+(\mathscr N, K)$ with itself $(h -1)$ times, and hence equal to $Z_+(\mathscr N, K).$
$\blacksquare$

\subsection{PY-TIK properties as PL complex balanced systems}

Am important property that can be inferred from the previous Section's results is that PY-TIK systems belong to the set of PY-RDK systems with zero kinetic deficiency:

\begin{proposition}
If $(\mathscr N, K)$ is a PY-TIK system, then $\widetilde\delta=0$.
\end{proposition}
\textit{Proof.}
According to the ZKRD Theorem, $(\mathscr N, K)$  is unconditionally complex balanced, i.e. $Z_+(\mathscr N, K)\neq\emptyset$ for any set of rate constants. According to Proposition \ref{prop:mure}, for PL complex balanced systems, this is equivalent to zero kinetic deficiency property.  
$\blacksquare$

\begin{remark}
An alternative proof is provided by using the representation of the PY-TIK system with PL-TIK systems, and using their Unconditional Complex Balancing to conclude that their kinetic deficiencies are all 0. The sum formula for kinetic deficiencies then shows the claim.
\end{remark}

In the case of power law systems $(h = 1)$, weakly reversible PL-RDK systems with zero kinetic deficiency display special properties. For example, Johnston \cite{JOHN2015} showed that complex balanced equilibria of systems with both network and kinetic deficiencies equal to zero are “translations” of toric steady states of mass action systems. We expect similar special properties in the poly-PL case, too.
In the following Sections, properties of PL complex balanced systems regarding linear stability and concentration robustness will be derived, which naturally also hold for PY-TIK systems.

\section{Linear stability of complex balanced PY-RDK systems} \label{sec:stable}

This section utilizes the STAR-MSC transform of a weakly reversible PY-RDK system in order to expand  the results of Boros et al. \cite{BMR2020} on linear stability to a subset of $Z_{+} (\mathscr{N},K)$ and a subset of PY-RDK system. We first recall some notions of square matrix stability  reviewed in \cite{BMR2020}. 
 
\begin{definition}
Let $S$ be a linear subspace. A square matrix $A$ with $\text{Im } A \subseteq S$ is said to be \textbf{stable on} $\bm{S}$ if all eigenvalues of the linear map $A|_S : S \rightarrow S$ have negative real part. If all eigenvalues of this linear map have non-positive real part, then the square matrix $A$ is \textbf{semistable on} $\bm{S}$.
\end{definition}

Denote the set of diagonal matrices with positive diagonal by $\mathcal{D}_+ \subseteq \mathbb{R}^{m \times m}$. 
\begin{definition}
Let $S$ be a linear subspace. A square matrix $A \in \mathbb{R}^{m \times m}$ with $\text{Im } A \subseteq S$  is  
\begin{enumerate}
\item \textbf{diagonally stable on } $\bm{S}$ (respectively, \textbf{diagonally semi-stable on }  $\bm{S}$) if there exists $P \in \mathcal{D}_+$ such that $PA+A^\top P <0$ on $S$ (respectively, $PA+A^\top P \leq 0$ on $S$)
\item  $\bm{D}$\textbf{-stable on} $\bm{S}$  (respectively, $\bm{D}$\textbf{-semistable on} $\bm{S}$) if $AD$ is stable on $S$ (respectively, semistable on $S$) for all $D \in \mathcal{D}_+$. 
\item \textbf{diagonally $\bm{D}$-stable on $\bm{S}$} (respectively, diagonally $D$-semistable on $S$) if, for all $D \in \mathcal{D}_+$, there exists $P \in \mathcal{D}_+$ such that $PAD + DA^\top P <0$ on $S$ (respectively, $PAD + DA^\top P \leq 0$ on $S$.
\end{enumerate}

\end{definition}

Recall that for a differentiable species formation rate function $f$, the \textbf{Jacobian matrix} of $f$ at $c$, denoted by $J(c)$,  is the $m \times m$ matrix where $
[J(c)]_{ji}= \dfrac{\partial f_j(c)}{\partial c_i}.
$

\begin{definition}
An equilibrium $c^*$ is \textbf{linearly stable} in its stoichiometric class $c^* + S$ if the Jacobian matrix $J(c^*)$ is stable on $S$. We say that $c^*$ is \textbf{diagonally \textit{D}-stable}, \textbf{diagonally stable}, or \textbf{\textit{D}-stable} in $c^* + S$ if $J(c^*)$ is diagonally \textit{D}-stable, diagonally stable, or \textit{D}-stable on $S$, respectively.
\end{definition}

One of the main results in \cite{BMR2020} was that linear stability of a complex balanced equilibrium of a GMAS system implies uniqueness in a stoichiometric class. The analogous result we have so far is the following:
\begin{theorem}
Let $(\mathscr{N},K)$ be a weakly reversible PY-RDK system.
\begin{enumerate}
\item[i.] If a PL-complex balanced equilibrium $c^*$ is linearly stable, then it is unique in its stoichiometric class.
\item[ii.] If $(\mathscr{N},K)$ is PL-complex balanced, then any linearly stable complex balanced equilibrium is unique in its stoichiometric class.
\end{enumerate}
\end{theorem}

\textit{Proof.} To prove (i), note that since $c^*$ is PL-complex balanced, it is a complex balanced equilibrium of $(\mathscr{N}^*,K^*)$. Since the latter is dynamically equivalent to $(\mathscr{N}, K)$, $c^*$ is a linearly stable equilibrium for $(\mathscr{N}^*,K^*)$. According to the GMAS result (i.e. Theorem 10 of \cite{BMR2020}), it is unique in its stoichiometric class under $\mathscr{N}^*$. However, since $S^* = S$, it is also unique in the stoichiometric class of $\mathscr{N}$.  Statement (ii) follows from the equation $Z_{+,\text{PL}}(\mathscr{N},K) = Z_+(\mathscr{N},K)$. $\blacksquare$

Finally, we briefly state the analogues of Theorems 11 and 13 from \cite{BMR2020} for PL-complex balanced PY-RDK systems on cyclic and weakly reversible networks respectively.
\begin{theorem}
Let $(\mathscr{N},K)$ be a PL-complex balanced PY-RDK system on a cyclic network. Let $A = YA_{k=1}\Psi_{K}$  ,where $\Psi_K$ is the factor map of $K$ and $A_{k=1}$ is the Laplacian map with all rate constants equal to 1. Then we have the following relationships:

$$
\small
\begin{array}
[c]{ccc}
\text{For all rate constants,} & & \text{For all rate constants,} \\
\text{complex balanced equilibria} & & \text{complex balanced equilibria} \\
\text{are diagonally stable} & \Rightarrow & \text{are linearly stable} \\
\text{in their stoichiometric classes.} & &\text{in their stoichiometric classes.} \\
\Updownarrow & & \Updownarrow \\
A \text{ is diagonally D-stable on } S. & \Rightarrow & A \text{ is D-stable on } S. 
\end{array}
$$
\end{theorem}
\vspace{0.5cm}
\begin{theorem}
Let $(\mathscr{N},K)$ be a PL-complex balanced PY-RDK system on a weakly reversible network. For each cycle $C$, let $A^C = YA^C_{k=1} \Psi_K$ and $S^C$ be the corresponding stoichiometric subspace. Then we have the following relationships:
$$
\small
\begin{array}
[c]{ccc}
\text{For all rate constants,} & & \text{For all rate constants,} \\
\text{complex balanced equilibria} & & \text{complex balanced equilibria} \\
\text{are diagonally stable} & \Rightarrow & \text{are linearly stable} \\
\text{in their stoichiometric classes.} & &\text{in their stoichiometric classes.} \\
\Downarrow & & \Downarrow \\
\text{For all cycles } C, & & \text{For all cycles } C, \\
A^C \text{ is diagonally D-semistable on } S^C . & \Rightarrow & A^C \text{ is D-semistable on } S^C .
\end{array}
$$
\end{theorem}

\section{Concentration robustness in PYK systems} \label{sec:robust}

In this section, we derive new results about concentration robustness in poly-PL systems after reviewing and consolidating previous ones. We begin with balanced concentration robustness (BCR), where the set of complex balanced equilibria is best understood and a simple computation procedure is available. We then consider sufficient conditions for absolute concentration robustness in PYK systems .

\subsection{Concentration robustness in poly-PL systems of CLP type} \label{sec:CRCLP}

We recall the general definition of concentration robustness in a species $X$ over a set of positive equilibria $Y$.

\begin{definition}
A kinetic system $(\mathscr{N}, K)$ has concentration robustness in a species $X$ over a subset $Y$ of positive equilibria if the concentration value for $X$ is invariant over all equilibria in $Y$. 
\end{definition}

Most studied sets have been absolute concentration robustness (ACR) where $Y = E_+(\mathscr{N},K)$ and balanced concentration robustness (BCR) when $Y=Z_+(\mathscr{N},K)$.

In the following, we first study balanced concentration robustness in PYK systems of CLP type, i.e. those with $Z_+(\mathscr{N},K) = \{ x \in \mathbb R^\mathscr S_> |\log x - \log x^*| \in P_Z^\perp\}$, with $x^*$ a given complex balanced equilibrium, $P_Z$ the system's flux subspace and $P_Z^\perp$ its parameter subspace. Lao et al derived the following ``Species Hyperplane Criterion for BCR'':

\begin{theorem}[Theorem 3.12, \cite{LLMM2021}]
Let $(\mathscr{N},K)$ be a kinetic system, $m_{\text{ACR}}$ and $m_{\text{BCR}}$ be the number of species with ACR and BCR, respectively. If $(\mathscr{N},K)$ is a weakly reversible CLP system then it has BCR in a species $X$ if and only if its parameter subspace $P^\perp_Z$ is a subspace of the species hyperplane $\{x\in\mathbb R^{\mathscr S}|x_X=0\}$. Furthermore, $m_{\text{ACR}}\leq m_{\text{BCR}}\leq s$.
\end{theorem}

The criterion above leads to a simple procedure for determining BCE in a weakly reversible CLP system. If $P^\perp_Z=0$, then the system has a unique complex balanced equilibrium, which is trivially equivalent to BCR in all species. Otherwise, one constructs a basis of $P^\perp_Z$, and  and BCR holds in $X$ if and only if all basis vectors have a zero coordinate in $X.$

In Proposition \ref{prop:mure}, we showed that every PL-complex balanced PY-RDK system is a CLP system with $P_Z = \tilde{S}_1 +\dots+\tilde{S}_h$.  Hence, balanced concentration robustness in such PYK systems can be determined by the results reviewed above.

Proposition 6.4 of \cite{LLMM2021} shows that the samd   conclusion as in Proposition (\ref{prop:mure}) holds even if the CLP summands are not necessarily PL-RDK. A class pf such PLK systems is given by the following result of \cite{FOME2021}:

\begin{theorem}
Let $(\mathscr{N},K)$ be a weakly reversible power law kinetic system with a complex balanced PL-RDK decomposition $\mathscr{D}: \mathscr{N}=\mathscr N_1\cup\cdots\cup\mathscr N_k$ with $P_{Z,i}=\tilde{S}_i$. If $\mathscr D$ is incidence independent and the induced covering $\tilde{\mathscr D}$ is independent, then $(\mathscr{N},K)$ is a weakly reversible CLP system with $P_Z=\sum \tilde{S}_i$.
\end{theorem}

A PL-NDK system with such PL-RDK decomposition is a subnetwork of Schmitz's earth pre-industrial carbon cycle model \cite{FMRL2019}. Such summands can be combined with PL-RDK systems to yield PL-complex balanced PY-NDK systems of CLP type.

If a PL-complex balanced PYK system is absolutely complex balanced, then ACR in a species if and only if BCR in a species. The case of systems with positive non-complex balanced equilibria is discussed in the next section.

\subsection{Absolute concentration robustness in poly-PL systems} \label{sec:revrobust}

We first introduce the analogue of PL-complex balanced PYK systems for the study of absolute concentration robustness.

\begin{definition}
A poly-PL kinetic system $ (\mathscr{N},K)$ is \textbf{PL-equilibrated (PLE)} if $E_{+} (\mathscr{N},K) =E_{+,PL} (\mathscr{N},K)$. The sets of PL-equilibrated PYK and PY-RDK systems are denoted by $\text{PYK}_\text{PLE}$ and $\text{PY-RDK}_\text{PLE}$.
\end{definition}

Our aim is to derive sufficient conditions for ACR in species of PL-equilibrated systems. These will be based on “low deficiency ACR building blocks” similar to those for PLK systems, which we now briefly recall.

\subsubsection{A brief review of ACR in low deficiency PL-RDK systems}

The concept of \textbf{absolute concentration robustness (ACR)} was first introduced by Shinar and Feinberg in their well-cited paper published in \textit{Science} \cite{SF2010}. ACR pertains to a phenomenon in which a species in a chemical kinetic system carries the same value for any positive steady state the network may admit regardless of initial conditions. The work of Shinar and Feinberg is influential largely because they established simple yet sufficient criteria for a mass action system to exhibit ACR. In \cite{FLRM2020}, this result is slightly modified to come up with an analogous theorem for PL-RDK systems.

Fortun and Mendoza \cite{FM2020} further investigated ACR in power law kinetic systems and derived novel results that guarantee ACR for some classes of PLK systems. For these PLK systems, the key property for ACR in a species $X$ is the presence of an SF-reaction pair. A pair of reactions in a PLK system is called a \textbf{Shinar-Feinberg pair} (or \textbf{SF-pair)} in a species $X$ if their kinetic order vectors differ only in $X$. A subnetwork of the PLK system is of \textbf{SF-type} if  it contains an SF-pair in $X$.

An SF-pair is called linked if both reactions lie in a linkage class. It is called non-terminal if both reactant complexes do not lie in terminal strong linkage classes.
The following sufficient condition for deficiency zero PL-RDK systems was derived in \cite{FLRM2020}:

\begin{theorem}[Th. 6, \cite{FM2020}] \label{th:acrdz}
Let $(\mathscr{N},K)$ be a deficiency zero PL-RDK with a positive equilibrium. If the system has a linked Sf-pair in a species $X$, then it has ACR in $X$.
\end{theorem}

We have the following sufficient condition for ACR in larger and higher deficiency PLK systems \cite{LLMM2021}:

\begin{proposition}[Prop. 5.3, \cite{LLMM2021}]
Let $(\mathscr{N},K)$ be a power law system with a positive equilibrium and an independent decomposition $\mathscr{N}=\mathscr N_1\cup\cdots\cup\mathscr N_k$. If there is a subnetwork $(\mathscr{N}_i,K_i)$ of deficiency $\delta_i$ with SF-pair in species $X$ such that
\begin{enumerate}
    \item [(i)] $\delta_i=0$ and $(\mathscr{N}_i,K_i)$ is a weakly reversible PL-RDK system with the SF-pair in a linkage class; or
    \item [(ii)] $\delta_i=0$, $(\mathscr{N}_i,K_i)$ is a PL-RDK system, and the SF-pairs's reactant complexes are nonterminal;
\end{enumerate}
then $(\mathscr{N},K)$ has ACR in $X$.
\end{proposition}

 \begin{theorem}[Shinar-Feinberg Theorem on ACR for PL-RDK systems, \cite{FLRM2020}] \label{th:SFTACR} Let $\mathscr{N}=(\mathscr{S,C,R})$ be a deficiency-one CRN and suppose that $(\mathscr{N},K)$ is a PL-RDK system which admits a positive equilibrium.  If $y, y' \in \mathscr{C}$ are nonterminal complexes whose kinetic order vectors  differ only in species $X$, then the system has ACR in $X$. \end{theorem}

\subsubsection{Low deficiency ACR building blocks for PYK systems}

The key step is the following extension of the SF-pair concept from PLK to PYK:

\begin{definition}\label{def:sf-pair}
Let $\{ \mathcal K_j \}$ the canonical PL-representation of the PYK system $(\mathscr{N}, K)$ and $R_i, R_{i'} \in \mathscr{R}$.  The pair of reactions $\{ R_i, R_{i'} \}$ is an \textbf{SF-pair in the species} $\bm{X}$ if there is at least one $\mathcal K_j$ where the rows differ in $X$. A PYK system with an SF-pair in $X$ is said to be of \textbf{Shinar-Feinberg type} (SF-type) in $X$. 
\end{definition}

An SF-pair is called linked if both reactions lie in a linkage class. It is called non-terminal if both reactant complexes do not lie in terminal strong linkage classes.

We derive the deficiency one building block first:

\begin{theorem}[SFACR Theorem for $\text{PY-RDK}_\text{PLE}$]\label{SFACR-PLE}
Let $\{ \mathcal K_j \}$ be the canonical PL-representation of a deficiency one $(\mathscr{N}, K)$ with $K \in \text{PY-RDK}_\text{PLE}$. If $E_+(\mathscr{N},K) \neq \emptyset$ and a non-terminal SF-pair in a species $X$ exists, then the system has ACR in $X$.
\end{theorem}

\textit{Proof.}
We first note that for any kinetic system $(\mathscr{N}, K)$ and a decomposition $\mathscr{N}=\mathscr N_1\cup\cdots\cup\mathscr N_k$ with $\emptyset\neq E_+(\mathscr{N}, K) = \cap E_+(\mathscr{N}_i,K_i)$, if $X$ has ACR in $(\mathscr{N}_i, K_i)$, then it has ACR in $(\mathscr{N}, K)$ too. Let $(\mathscr{N}^*, K^*)$ be the STAR-MSC transform of $(\mathscr{N}, K)$. We have $E_+(\mathscr{N}, \mathcal K_j)=E_+(\mathscr{N}_j^*, K^*)$, hence, due to dynamic equivalence, $E_+(\mathscr{N}, K)=E_+(\mathscr{N}^*, K^*)$, and PL-equilibrated implies 
$E_+(\mathscr{N}^*, K^*)=\cap E_+(\mathscr{N}_j^*, K^*)$. There is at least one $j$ such that a non-terminal SF-pair for $(\mathscr{N}_j^*, K^*)$ exists for the deficiency 1 network $\mathscr{N}_j^*$ which has a positive equilibrium for the PL-RDK kinetics $K^*$. Hence the subnetwork $\mathscr{N}^*_j$ has ACR in $X$. It follows that the whole network has ACR in $X.$
$\blacksquare$

An analogous proof can be given for the deficiency zero building block:

\begin{theorem}[ACR Theorem for zero deficiency $\text{PY-RDK}_\text{PLE}$]\label{ACR-PLE}
Let $\{ \mathcal K_j \}$ be the canonical PL-representation of a weakly reversible deficiency zero $(\mathscr{N}, K)$ with $K \in \text{PY-RDK}_\text{PLE}$. If $E_+(\mathscr{N},K) \neq \emptyset$ and a non-terminal linked SF-pair in a species $X$ exists, then the system has ACR in $X$.
\end{theorem}

\subsubsection{ACR in larger and higher deficiency PYK systems}

We can now provide a sufficient condition for ACR in larger and higher deficiency PYK systems: 

\begin{proposition}
Let $(\mathscr{N},K)$ be a poly PL- system with a positive equilibrium and an independent decomposition $\mathscr{N}=\mathscr N_1\cup\cdots\cup\mathscr N_k$. If there is a subnetwork $(\mathscr{N}_i,K_i)$ of deficiency $\delta_i$ with SF-pair in species $X$ such that
\begin{enumerate}
    \item [(i)] $\delta_i=0$ and $(\mathscr{N}_i,K_i)$ is a weakly reversible $\text{PY-RDK}_\text{PLE}$ system with the SF-pair in a linkage class; or
    \item [(ii)] $\delta_i=0$, $(\mathscr{N}_i,K_i)$ is a $\text{PY-RDK}_\text{PLE}$ system, and the SF-pairs's reactant complexes are nonterminal;
\end{enumerate}
then $(\mathscr{N},K)$ has ACR in $X$.
\end{proposition}

\textit{Proof.}
Since the decomposition is independent, it follows from Feinberg's Decomposition Theorem that $E_+(\mathscr{N},K)=\cap E_+(\mathscr{N}_i,K_i)$. ACR in $X$ in the superset of equilibria $E_+(\mathscr{N}_i,K_i)$ of a building block implies the claim.
$\blacksquare$

\section{Summary and Outlook}
\label{sec:summary}
Using tools in chemical reaction network theory (CRNT), this paper assembles relevant results that tackle the dynamical properties of chemical reaction networks endowed with poly-PL kinetics. Poly-PL kinetic systems consist of nonnegative linear combinations of power law functions. The main approach used in this paper is the construction of power law representations of poly-PL systems in order to apply, expand or mimic existing results  in CRNT involving power law kinetic systems. This contribution offers two primary methods in rewriting a poly-PL system into its power law representation: the canonical PL-representation and the STAR-MSC representation. Based from these representations, the following main results are reiterated:
\begin{enumerate}
\item For a weakly reversible PY-RDK system, there are rate constants for which the system is complex balanced. In fact, if its canonical PL-representation consists of PL-RDK systems with zero kinetic deficiency, then the system is complex balanced for any set of rate of constants. The parametrization of a subset of the set of complex balanced equilibria is provided in Proposition \ref{prop:gmak2}. Further, a necessary and sufficient condition for the existence of more than one distinct complex balanced PL-equilibria in a stoichiometric class of a weakly reversible PY-RDK system is established.  All these aforementioned results form partial extensions of the GMAS results of M\"{u}ller and Regensburger \cite{MURE2012,MURE2014} to a subset of poly-PL kinetics that are reactant-determined (PY-RDK). Full analogue results for PY-RDK systems with PL-equilibria and PL-complex balanced equilibria sets are collated and presented in Proposition \ref{prop:mure}. 

\item A class of poly-PL systems called PY-TIK systems that are weakly reversible were shown to exhibit PL-complex balancing property. They are precisely the poly-PL systems with zero kinetic reactant deficiency. 

\item STAR-MSC transform of a weakly reversible poly-PL kinetic system has led us to the extension of the results of Boros et al. \cite{BMR2020} on linear stability to a subset of complex balanced equilibria and a subset of PY-RDK system. 
 
\item The notions of absolute concentration robustness (ACR) and balanced concentration robustness (BCR) are also extended to PL-equilibrated and PL-complex balanced poly-PL kinetic systems, respectively. The results are expansions of the work of Fortun and Mendoza \cite{FM2020} on deficiency zero ACR,  and the detection of ACR and BCR in larger or higher deficiency networks.

 \end{enumerate}

As future perspective, one may look at the extension of the multistationarity algorithm for power law kinetic systems \cite{HERNZ2020} to poly-PL systems. The results provided in this paper can be further enlarged to include Hill-type kinetic systems through the transformation introduced in \cite{HM2020}.




\vspace{0.5cm}
\noindent \textbf{Acknowledgement.}
We thank Dr. Noel T. Fortun for his contributions to an earlier version of the manuscript. ECJ's work was funded by the UP System Enhanced Creative Work and Research Grant (ECWRG 2020-1-7-R).

\newpage
\appendix
\section{Complex Balancing of RID Kinetic Systems} \label{append:A}

We derive the Conditional Complex Balancing (CCB) property for any poly-PL kinetics on a weakly reversible network from a general result on RID kinetics.

\begin{definition}
A \textbf{rate constant-interaction map decomposable kinetics} (RIDK) is a kinetics, such that for each reaction $R_j$, the coordinate function $K_j: \Omega \rightarrow \mathbb{R}$ can be written in the form $K_j(x)=k_j I_{K,j} (x)$, with $k_j \in \mathbb{R}_{>0}$ (called rate constant) and $\Omega \subset \mathbb{R}^m$. We call the map $I_K : \Omega \rightarrow \mathbb{R}^r$ defined by $I_{K,j}$ as the \textbf{interaction map}.
\end{definition}

A poly-PL kinetics $K$ is an RIDK, the interaction $I_K$ being the poly-PL function. By definition, an RIDK is called complex factorizable if at each branching point of the network, any branching reaction has the same interaction \cite{NEML2019}.

We recall three known facts in a lemma:

\begin{lemma}\label{lem:ccb}
Consider the chemical kinetic system $(\mathscr{N},K)$.
\begin{enumerate}
\item[i.] A network $\mathscr{N}$ is weakly reversible if and only if $\text{Ker } I_a$ contains a positive vector.
\item[ii.] A weakly reversible network $\mathscr{N}$ is positive dependent.
\item[iii.] For any RIDK on a positively dependent network $\mathscr{N}$, there are rate constants so that $(\mathscr{N},K)$ has a positive equilibrium.
\end{enumerate}
\end{lemma}

\noindent The proofs of (i) and (ii) can be found in \cite{DELE2009} and \cite{FEIN1987}, respectively. The proof in \cite{FEIN1987} is readily adapted to show (iii).

\begin{theorem}[CCB for RIDK on weakly reversible networks]
Let $K$ be an RID kinetics on a weakly reversible network $\mathscr{N}$.  Then there exist rate constants such that $(\mathscr{N}, K)$ is complex balanced, i.e. $Z_+(\mathscr{N},K) \neq \emptyset$.
\end{theorem}

\textit{Proof.}
Since a weakly reversible network is positive dependent, according to Lemma \ref{lem:ccb} (ii) and (iii), there are rate constants $k^*_q$ such that $(\mathscr{N}, K)$ has a positive equilibrium $x^*$. Let $I_K$ be the interaction function of $K$. Furthermore, according to Lemma \ref{lem:ccb} (i), weak reversibility implies that there is a positive vector $b \in \text{Ker } I_a$.  Let $C = \text{diag } (b_1/ I_{K,1}(x^*),\dots, b_r/ I_{K,r}(x^*))$ and $k_q' = (Ck^*)_q$. Then $I_{a,q} (k_q'I_{K,q}(x^*)) = I_{a,q}(b_q/ I_{K,q}(x^*))k_q^* I_{K,q}(x^*) = k_q^*I_a(b) = 0$. Hence, for rate constants $ k_q'$, $x^*$ is a complex balanced equilibrium, and the claim is established.
$\blacksquare$

\section{Nomenclature} \label{append:B}
We list some of the symbols and acronyms used in the paper.

\vspace{0.5cm}
\noindent \textbf{List of Symbols}
\begin{table}[ht!]
\begin{tabular}{ll}
$\delta$	& 	Deficiency	\\
$K$	& 	Kinetics of a CRN	\\
$\Psi_K$	& 	Factor map of a kinetics	\\
$\text{\textsf{sum}} \bm{\left\langle B_j \right\rangle}$	& 	Image span sums	\\
$I_a$	& 	Incidence map/matrix	\\
$J$	& 	Jacobian matrix	\\
$\widetilde{\delta}$	& 	Kinetic deficiency	\\
$F$	& 	Kinetic order matrix	\\
$A_k$	& 	$k$-Laplacian matrix 	\\
$Y$	& 	Map of complexes	\\
$Z_+(\mathscr{N},K)$	& 	Set of complex balanced equilibria	\\
$E_{+,PL}(\mathscr{N},K)$	& 	Set of Plcomplex balanced	\\
$E_{+,PL}(\mathscr{N},K)$	& 	Set of PL-equilibria	\\
$E_+(\mathscr{N},K)$	& 	Set of positive equilibria	\\
$N$	& 	Stoichiometric matrix	\\
\end{tabular}
\end{table}

\newpage
\noindent \textbf{Abbreviations}
\begin{table}[ht!]
\begin{tabular}{ll}
ACR	&	Absolute concentration robustness	\\
BCR	&	Balanced concentration robustness	\\
CKS	&	Chemical kinetic system	\\
CRN	&	Chemical reaction network	\\
CF	&	Complex factorizable	\\
CCB	&	Conditional complex balancing	\\
GMAS	&	Generalized mass action system	\\
MAK	&	Mass action kinetic	\\
NF	&	Non-complex factorizable	\\
NDK 	&	Non-reactant determined kinetics	\\
PLC	&	PL-complex balanced	\\
PLE	&	PL-equilibrated	\\
PY	&	Poly-PL	\\
PYK	&	Poly-PL kinetic system	\\
PL	&	Power law	\\
RID	&	Rate constant-interaction map decomposable	\\
RDK	&	Reactant-determined kinetics	\\
SF	&	Shinar-Feinberg	\\
STAR-MSC	&	$S$-invariant termwise addition of reactions via maximal stoichiometric coefficients	\\
UCB	&	Unconditional complex balancing	\\
\end{tabular}
\end{table}

\end{document}